\documentclass[12pt,a4paper,english]{article}
\usepackage{babel,amsmath,amssymb,amsthm}

\newtheorem*{claim*}{Claim}
\newtheorem*{definition*}{Definition}
\newtheorem*{lemma*}{Lemma}

\newcommand{\EHR}{\mathrm{EHR}}

\title{Bounded quantifier depth spectrum for random uniform hypegraphs}
\author{}
\date{}
\begin{document}
\maketitle

\begin{center}
{\large S.N. Popova~\footnote{Moscow Institute of Physics and Technology, Laboratory of Advanced Combinatorics and Network Applications; National Research University Higher School of Economics.}}
\end{center}

The notion of spectrum for first-order properties introduced by J. Spencer for Erdős–Rényi random graph 
is considered in relation to random uniform hypergraphs. In this work we study the set of limit points 
of the spectrum for first-order formulae with bounded quantifier depth and obtain bounds for its maximum value. 
Moreover, we prove zero-one $k$-laws for the random uniform hypergraph and improve the bounds
for the maximum value of the spectrum for first-order formulae with bounded quantifier depth. We obtain that the maximum value of the spectrum belongs to some two-element set. 


\section{Introduction}
\indent 

Asymptotic behavior of first-order properties probabilities of the Erdős–Rényi random graph
$G(n, p)$ has been widely investigated (see \cite{shelah_spencer}--\cite{raigor}).
In \cite{inf_spectra} the notion of spectrum of first-order properties was introduced and it was proved that there exists a first-order property with infinite spectrum.
In this work we consider spectrum of first-order properties in relation to random uniform hypergraphs.
Let us define the \textit{random $s$-uniform hypergraph} $G^s(n, p)$. Consider the set $\Omega_n = \{G = (V_n, E)\}$ of all $s$-uniform hypergraphs ($s$-hypergraphs) 
with the set of vertices $V_n = \{1, 2, \ldots, n\}$. The random hypergraph $G^s(n, p)$ is a random element with probability distribution
$$\Pr[G^s(n, p) = G] = p^{|E|} (1 - p)^{{n \choose s} - |E|}.$$
Note that for $s = 2$ in this definition we obtain the Erdős–Rényi random graph $G(n, p)$.
Let us denote the event ``$G^s(n, p)$ has a property $L$'' by $G^s(n, p) \models L$.
\textit{First-order properties} of $s$-uniform hypergraphs are defined by first-order formulae (see \cite{raigor}, \cite{ebbing}) which are built of predicate symbols $N, =$, logical connectivities, variables and quantifiers $\forall, \exists$. The $s$-ary predicate symbol $N$ expresses the property of $s$ vertices to constitute an edge.  
Let us recall that the \textit{quantifier depth} (see \cite{raigor}, \cite{ebbing}) of a first-order formula is the maximum number of nested quantifiers.
Let $\mathcal L_k$ denote the set of all $s$-hypergraph properties which can be expressed by first-order formulae with quantifier depth at most $k$. Moreover, let $\mathcal L = \cup_{k \in \mathbb N} \mathcal L_k$ be the set of all first-order $s$-hypergraph properties. 

The random hypergraph $G^s(n, p)$ is said to \textit{obey the zero-one law} if for any first-order property $L \in \mathcal L$ the probability $\Pr[G^s(n, p) \models L]$ tends either to 0 or to 1 as $n \to \infty$. We say that the random hypergraph $G^s(n, p)$ \textit{obeys the zero-one $k$-law} if for any first-order property $L \in \mathcal L_k$ the probability $\Pr[G^s(n, p) \models L]$ tends either to 0 or to 1 as $n \to \infty$.

Let us define the \textit{spectrum} $S(L)$ for any first order property $L \in \mathcal L$. 
$S(L)$ is the set of all $ \alpha \in (0, s - 1)$ 
(we take the interval $(0, s - 1)$ because the structure of $G^s(n, n^{-\alpha})$ for $\alpha > s - 1$  is considerably simpler than that for $\alpha < s- 1$ and the study of $G^s(n, p)$ with $p(n) \gg n^{- s + 1}$ is a subject of a separate investigation) which do not satisfy the
following property: $\lim_{n \to \infty} \Pr[G^s(n, p) \models L]$ exists and is either $0$ or $1$.
Denote the union of $S(L)$ over all $L\in \mathcal L_k$ by $S_k$.
In other words, $S_k$ is the set of all $\alpha \in (0, s - 1)$ for which the random hypergraph $G^s(n, n^{-\alpha})$ does not obey the zero-one $k$-law. 

Let $(S_k)'$ be the set of limit points in $S_k$.

In \cite{shelah_spencer} S. Shelah and J. Spencer proved that when $\alpha$ is an irrational number and $p(n) = n^{-\alpha + o(1)}$ then the random graph $G(n, p)$ obeys the zero-one law. If $\alpha \in (0, 1) \cap \mathbb Q$, where $\mathbb Q$ is the set of all rational numbers, then $G(n, n^{-\alpha})$ does not obey the zero-one law. If $\alpha > 1$, then $G(n, n^{-\alpha})$ obeys the zero-one law if and only if $\alpha \notin \{1 + \frac{1}{k}: k \in \mathbb N\}$. 

In \cite{zhuk1}-\cite{zhuk4} zero-one $k$-laws for the random graph $G(n, p)$ were investigated. 
It was proved that if $\alpha$ is close to $0$ or close to $1$, then $G(n, n^{-\alpha})$ obeys the zero-one $k$-law. The spectrum $S_k$ for the random graph $G(n, p)$ was studied and it was shown that
$\min S_k = \frac{1}{k-2}$ and $\max S_k = 1 - \frac{1}{2^{k-2}}$. 
Moreover, $\min (S_k)' \le \frac{1}{k-11}$ and $\max (S_k)' \ge 1 - \frac{1}{2^{k-5}}$.

In \cite{zero_one} the result from \cite{shelah_spencer} was extended to the case of random uniform hypergraph $G^s(n, p)$. The random hypergraph $G^s(n, n^{-\alpha})$ obeys the zero-one law if and only if $\alpha \in (0, s - 1) \setminus \mathbb Q$ or $\alpha \in (s-1, \infty) \setminus \{s - 1 + \frac{1}{k}: k \in \mathbb N\}$.

In our previous papers we studied the spectrum $S_k$ for the random hypergraph $G^s(n, p)$ and examined the question for which $k$ the set $S_k$ is infinite (see \cite{infspectrum}). We also estimated the minimum and maximum values in $S_k$ and proved zero-one $k$-laws for the random hypergraph $G^s(n, p)$ (see \cite{spectrum}, \cite{zero_one}, \cite{limitspectrum}). 

We showed that there exists an interval with the left endpoint equals to $0$, such that for all $\alpha$ from this interval the random hypergraph $G^{s}(n,n^{-\alpha})$ obeys the zero-one $k$-law.

\textbf{Theorem 1.} (\cite{spectrum}) \textit{
Let $s \ge 3$, $k \ge s + 1$, $\alpha > 0$ and 
$$\frac{1}{\alpha} > {{k - 1} \choose {s - 1}} - 1 - \frac{s - 1}{k - 1} + \frac{2\left(1 + \frac{s - 1}{k - 1}\right)}{{{k - 1} \choose {s - 1}} + 2}.$$
Then the random hypergraph $G^s(n, n^{-\alpha})$ obeys the zero-one $k$-law.
}\\

We also proved that near the right endpoint of this interval there is a point $\alpha$ for which $G^{s}(n,n^{-\alpha})$ does not obey the zero-one $k$-law.

\textbf{Theorem 2.}  (\cite{spectrum}) \textit{
Let $s \ge 3$, $k \ge s + 2$. Then there exists $\alpha > 0$ such that 
$$\frac{1}{\alpha} > {{k - 1} \choose {s - 1}} - 1 - \frac{s - 1}{k - 1} - \frac{2}{{{k - 1} \choose {s - 1}}}$$
and $G^s(n, n^{-\alpha})$ does not obey the zero-one $k$-law.
}\\

These theorems provide a tight bound for the minimal point in $S_k$.
The question what is the exact value of the minimal point for which the zero-one $k$-law does not hold remains open.
Theorems 1 and 2 imply that this value asymptotically equals
$1/\left({{k - 1} \choose {s - 1}} - 1 - \frac{s - 1}{k - 1} + O\left(\frac{1}{{{k - 1} \choose {s - 1}}}\right) \right)$ as $k \to \infty$. 

Furthermore, we examined the zero-one $k$-law for the left neighborhood of $s-1$.

\textbf{Theorem 3}. (\cite{spectrum}) \textit{
Let $\alpha \in \left(s - 1 - \frac{1}{2^{k - s + 1}}, s - 1\right) \setminus \mathcal Q_k$, where
$\mathcal Q_k = \{s - 1 - \frac{1}{2^{k - s + 1} + \frac{a}{b}} | a, b \in \mathbb N, a \le 2^{k - s + 1}\}$.
Then the random hypergraph $G^s(n, n^{-\alpha})$ obeys the zero-one $k$-law. 
}\\

\textbf{Theorem 4}.  (\cite{spectrum}) \textit{
Let $s \ge 3$, $k \ge s + 4$, $\alpha = s - 1 - \frac{1}{2^{k - s + 1} + a}$, where $a \in \mathbb N$, $a \le 2^{k - s - 2} + 2^{k - s - 3} + 1$.
Then the random hypergraph $G^s(n, n^{-\alpha})$ does not obey the zero-one $k$-law.
}\\

From Theorems 3 and 4 we obtain that the maximal point of $S_k$ lies between $s - 1 - \frac{1}{2^{k-s+1} + 2^{k-s-2} + 2^{k-s-3} + 1}$ and $s - 1 - \frac{1}{2^{k-s+2}}$.\\

We also estimated the minimum value in $(S_k)'$.

\textbf{Theorem 5.} (\cite{limitspectrum}) \textit{
There exists a constant $C$ such that for any $k \ge s + C$ we have $\min (S_k)' \le \frac{1}{{{k - 11} \choose {s - 1}}}$.
}\\

Taking into account the bound from Theorem 1 we deduce that Theorem 5 gives an asymptotically tight bound for $\min (S_k)'$.

\textbf{Corollary 1.} \textit{
We have $\min (S_k)' \sim \frac{(s-1)!}{k^{s - 1}}$, as $k \to \infty$.
}\\

\section{New results}
\indent

We examine the set $(S_k)'$ of limit points in $S_k$ and estimate the maximum value in $(S_k)'$.

\textbf{Theorem 6.} \textit{
There exists a constant $C$ such that for any $k \ge s + C$ we have $\max (S_k)' \ge s - 1 - \frac{1}{2^{k - s - 4}}$.
}\\

From Theorems 3, 4 and 6 it follows that $- \log_2 (s - 1 - \max (S_k)') = k - s + O(1)$.

We prove the zero-one $k$-law for the random hypergraph $G^s(n, n^{-\alpha})$ where $\alpha$ belongs to some set of rational numbers from
$\left(s - 1 - \frac{1}{2^{k - s + 1}}, s - 1\right)$.

\textbf{Theorem 7.} \textit{Let $s \ge 3$, $k \ge s + 1$ and $\frac{a}{b}$ be an irreducible positive fraction. Denote $\nu = \max(1, 2^{k - s + 1} - b)$. Let $a \in \{\nu, \nu + 1, \ldots, 2^{k - s + 1}\}$, $\alpha = s - 1 - \frac{1}{2^{k - s + 1} + \frac{a}{b}}$.
Then the random hypergraph $G^s(n, n^{-\alpha})$ obeys the zero-one $k$-law. 
}\\

This theorem extends the class of rational fractions $\alpha$ for which the random hypergraph $G^s(n, n^{-\alpha})$ obeys the zero-one $k$-law and improves the bound obtained in Theorem 3. Theorem 7 implies that the maximal point of $S_k$ is not greater than $s - 1 - \frac{1}{2^{k-s+2} - 2}$. 
From Theorem 7 it follows that $S_k \cap \left(s - 1 - \frac{1}{2^{k - s + 1}}, s - 1\right)$ is finite, so $\max (S_k)' \le s - 1 - \frac{1}{2^{k - s + 1}}$. 

We also disprove the zero-one $k$-law for the random hypergraph $G^s(n, n^{-\alpha})$ where $\alpha$ belongs to the left neighborhood of $s-1$.

\textbf{Theorem 8}. \textit{Let $s \ge 3$, $k \ge s + 1$ and $\alpha = s - 1 - \frac{1}{2^{k - s + 1} + a}$, where $a \in \mathbb N$, $a \le 2^{k - s + 1} - 3$. Then the random hypergraph $G^s(n, n^{-\alpha})$ does not obey the zero-one $k$-law.
}\\

Theorem 8 improves the bound for the maximal point of $S_k$ which follows from Theorem 4.
Theorems 7 and 8 imply that $\max(S_k) \in \bigl\{s - 1 - \frac{1}{2^{k-s+2} - 3}, s - 1 - \frac{1}{2^{k-s+2} - 2}\bigr\}$.

\section{Auxiliary statements}
\indent

\subsection{Small subhypergraphs}
\indent

For an arbitrary $s$-hypergraph $G = (V, E)$, set $v(G) = |V|$, $e(G) = |E|$, $\rho(G) = \frac{e(G)}{v(G)}$, $\rho^{\max}(G) = \max_{H \subseteq G} \rho(H)$. $G$ is called \textit{strictly balanced} if the \textit{density} $\rho(G)$ of this graph is greater than the density of any its proper subhypergraph.
Denote the property of containing a copy of $G$ by $L_G$.

\textbf{Theorem 9} (\cite{vantsyan}). \textit{If $p \ll n^{-1/\rho^{\max}(G)}$, then
$\lim_{n \to \infty} \Pr[G^s(n, p) \models L_G] = 0$.
If $p \gg n^{-1/\rho^{\max}(G)}$, then
$\lim_{n \to \infty} \Pr[G^s(n, p) \models L_G] = 1$.}\\

In other words, the function $n^{-1/\rho^{\max}(G)}$ is a threshold for the property $L_G$.

Let $G_1, \dots, G_m$ be strictly balanced $s$-hypergraphs, $\rho(G_1) = \ldots = \rho(G_m) = \rho$. Let $a_i$ be the number of automorphisms of $G_i$.
Denote by $N_{G_i}$ the number of copies of $G_i$ in $G^{s}(n,p)$. The following theorem is a generalization of a classical result of Bollob\'{a}s (see, for example,~\cite{JLR}) to the case of $s$-uniform hypergraphs.

\textbf{Theorem 10}. \textit{
If $p=n^{-1/\rho}$, then
$$
(N_{G_1}, \dots, N_{G_m}) \xrightarrow{d}(P_1, \dots, P_m),
$$
where $P_i \sim \mathrm{Pois}\left(\frac{1}{a_i}\right)$ are independent Poisson random
variables.
}\\

\subsection{Extensions}
\indent

Consider arbitrary $s$-hypergraphs $G$ and $H$ such that $H \subset G$, $V(H) = \{x_1, \ldots, x_l\}$, $V(G) = \{x_1, \ldots, x_m\}$. Denote $v(G, H) = v(G) - v(H)$, $e(G, H) = e(G) - e(H)$,  $\rho(G, H) = \frac{e(G, H)}{v(G, H)}$,
$\rho^{\max}(G, H) = \max_{H \subset K \subseteq G} \rho(K, H)$. 
The pair $(G, H)$ is called \textit{strictly balanced} if $\rho(G, H) > \rho(K, H)$ for any graph $K$ such that
$H \subset K \subset G$.

Consider $s$-hypergraphs $\tilde H, \tilde G$, where $V(\tilde H) = \{\tilde x_1, \ldots, \tilde x_l\}$, $V(\tilde G) = \{\tilde x_1, \ldots, \tilde x_m\}$,
$\tilde H \subset \tilde G$. The hypergraph $\tilde G$ is called \textit{$(G, (x_1, \ldots, x_l))$-extension} of the tuple $(\tilde x_1, \ldots, \tilde x_l)$, if
$$
\{x_{i_1}, \ldots, x_{i_s}\} \in E(G) \setminus E(H) \Rightarrow \{\tilde x_{i_1}, \ldots, \tilde x_{i_s}\} \in E(\tilde G) \setminus E(\tilde H).
$$ 
If
$$
\{x_{i_1}, \ldots, x_{i_s}\} \in E(G) \setminus E(H) \Leftrightarrow \{\tilde x_{i_1}, \ldots, \tilde x_{i_s}\} \in E(\tilde G) \setminus E(\tilde H),
$$ 
we call $\tilde G$ a \textit{strict $(G, (x_1, \ldots, x_l))$-extension} of $(\tilde x_1, \ldots, \tilde x_l)$.

Let $\alpha > 0$. For any pair $(G, H)$, where $H \subset G$, set $f_\alpha(G, H) = v(G, H) - \alpha e(G, H)$. 
The pair $(G, H)$ is called \textit{$\alpha$-safe}, if
$f_\alpha(K, H) > 0$ for any $K$, $H \subset K \subseteq G$. If for any $K$ such that $H \subset K \subseteq G$ we have
$f_\alpha(G,K) < 0$, then the pair $(G, H)$ is called \textit{$\alpha$-rigid}. If for any $K$ such that $H \subset K \subset G$ we have
$f_\alpha(K, H) > 0$ and $f_\alpha(G, H) = 0$, then the pair $(G, H)$ is called \textit{$\alpha$-neutral}.

Let $\tilde H \subset \tilde G \subset \Gamma$, $T \subset K$ and $|V(T)| \le |V(\tilde G)|$. The pair $(\tilde G, \tilde H)$ is called \textit{$(K, T)$-maximal} in $\Gamma$, if for any subhypergraph $\tilde T \subset \tilde G$ with $|V(\tilde T)| = |V(T)|$ and $\tilde T \cap \tilde H \neq \tilde T$ there does not exist
a strict $(K, T)$-extension $\tilde K$ of $\tilde T$ in the hypergraph $\Gamma \setminus (\tilde G \setminus \tilde T)$ such that $E((\tilde K \cup \tilde G) \setminus \tilde T) \setminus (E(\tilde K \setminus \tilde T) \cup
E(\tilde G \setminus \tilde T)) = \varnothing$.
The hypergraph $\tilde G$ is called \textit{$(K, T)$-maximal} in $\Gamma$, if for any subhypergraph $\tilde T \subset \tilde G$ with $|V(\tilde T)| = |V(T)|$
there does not exist
a strict $(K, T)$-extension $\tilde K$ of $\tilde T$ in the hypergraph $\Gamma \setminus (\tilde G \setminus \tilde T)$ such that $E((\tilde K \cup \tilde G) \setminus \tilde T) \setminus (E(\tilde K \setminus \tilde T) \cup
E(\tilde G \setminus \tilde T)) = \varnothing$.

Let a pair $(G, H)$ be $\alpha$-safe and $\mathcal K_r$ be the set of all $\alpha$-rigid and $\alpha$-neutral pairs $(K, T)$, where
$|V(T)| \le |V(G)|$, $|V(K) \setminus V(T)| \le r$. Let $\tilde x_1, \ldots, \tilde x_l \in V_n$.
Define a random variable $N^{r}_{(G, H)}(\tilde x_1, \ldots, \tilde x_l)$ which assigns to each hypergraph $\mathcal G \in \Omega_n$ the number of all strict
$(G, H)$-extensions $\tilde G$ of the hypergraph $\tilde H = \mathcal G|_{\{\tilde x_1, \ldots, \tilde x_l\}}$ such that for any $(K, T) \in \mathcal K_r$ the pair $(\tilde G, \tilde H)$ is $(K, T)$-maximal in $\mathcal G$.

The following theorem is the generalization of Lemma 10.7.6 from \cite{alon_spencer} and Proposition 1 from \cite{zhuk_ext} to the case of $s$-uniform hypergraphs. The proof of this theorem is analogous to the proof for graphs.

\textbf{Theorem  11.}
\textit{
Let $p(n) = n^{-\alpha}$, $r \in \mathbb N$ and $(G, H)$ be $\alpha$-safe.
Then a.a.s. for every vertices $\tilde x_1, \ldots, \tilde x_l$ the following relation holds:
$$
N^{r}_{(G, H)}(\tilde x_1, \ldots, \tilde x_l) \sim {\sf E}[N^{r}_{(G, H)}(\tilde x_1, \ldots, \tilde x_l)] = \Theta(n^{f_\alpha(G, H)}).
$$
}\\

Denote by $L^{r}_{(G, H)}$ the following property: for any tuple $(\tilde x_1, \ldots, \tilde x_l)$ there exists a strict $(G, (x_1, \ldots, x_l))$-extension which is $(K, T)$-maximal for all $(K, T) \in \mathcal K_r$. 
Theorem 11 implies that if $(G, H)$ is $\alpha$-safe, then the random hypergraph $G^s(n, n^{-\alpha})$ satisfies $L^r_{(G, H)}$ a.a.s.

Denote by $\tilde L_{(G, H)}$ the property that there exists a copy of $H$ such that no copy of $G$ contains it.
Denote by $\tilde N_{(G, H)}$ the number of copies of $H$ in $G^s(n, p)$ which are not contained in any copy of $G$.

\textbf{Proposition 1} (\cite{infspectrum}). \textit{Let $p(n) = n^{-\alpha}$, $H$ be a strictly balanced $s$-hypergraph, $(G, H)$ be a strictly balanced pair, $\rho(H) = \rho(G, H) = 1/\alpha$.
Let $a_1$ be the number of automorphisms of $H$ which are extendable to some automorphism of $G$. Let $a_2$ be the number of automorphisms 
$\sigma: V(G) \to V(G)$ with $\sigma(x) = x$ for all $x \in V(H)$.
Then $$\tilde N_{(G, H)} \stackrel{d}{\longrightarrow} \mathrm{Pois}\left(\frac{1}{a(H)}\exp\left(-\frac{a(H)}{a_1 a_2}\right)\right).$$}\\

\subsection{Cyclic extensions}
\indent

We say that $(G,H)$ is a {\it cyclic $m$-extension} if $\rho^{\max}(G)<\frac{m}{m(s-1)-1}$, and $(G,H)$ fits one of the following patterns.

\begin{itemize}
\item $V(G,H) = \{y_1, \ldots, y_{k(s - 1)}, z_1, \ldots, z_l\}$, and there exists a vertex $x_1 \in H$ such that
\begin{multline*}
E(G,H) = \{\{x_1, y_1, \ldots, y_{s - 1}\}, \{y_{s-1}, \ldots, y_{2(s-1)}\}, \ldots, \{y_{(k - 1)(s - 1)}, \ldots, y_{k(s - 1)}\}, \\
\{y_{k(s - 1)}, z_1, \ldots, z_l, u_1, \ldots, u_{s - 1 - l}\}\},
\end{multline*}  
where $1 \le k \le m - 1$, $0 \le l < s - 1$ and $u_1, \ldots, u_{s - 1 - l} \in \{x_1, y_1, \ldots, y_{k(s - 1) - 1}\}$ are distinct. 
In such a situation, $(G, H)$ is called \textit{the first type} extension.

\item $V(G,H) = \{y_1, \ldots, y_{k(s - 1)}, z_1, \ldots, z_l\}$, and there exist distinct vertices $x_1,x_2 \in H$ such that
\begin{multline*}
E(G,H) = \{\{x_1, y_1, \ldots, y_{s - 1}\}, \{y_{s-1}, \ldots, y_{2(s-1)}\}, \ldots, \{y_{(k - 1)(s - 1)}, \ldots, y_{k(s - 1)}\}, \\
\{x_2, y_{k(s - 1)}, z_1, \ldots, z_l, u_1, \ldots, u_{s - 2 - l}\}\},
\end{multline*}
where $1 \le k \le m - 1$, $0 \le l \le s - 2$ and $u_1, \ldots, u_{s - 2 - l} \in \{y_1, \ldots, y_{k(s - 1) - 1}\}$ are distinct.
If $x_1 \neq x_2$, then $(G, H)$ is called \textit{the second type} extension.

\item $V(G,H) = \{y_1, \ldots, y_{s - l}\}$, $2 \leq l \leq s-1$, and there exist distinct vertices $x_1, \ldots, x_l \in H$, such that
$$
E(G,H) = \{x_1, \ldots, x_l, y_1, \ldots, y_{s - l}\}.
$$
In such a situation, $(G, H)$ is called \textit{the second type} extension.
\end{itemize}

For any $m \ge 1$, define a set of $s$-hypergraphs $\mathcal H_m$. The one-vertex hypergraph $(\{x\}, \varnothing)$ belongs to $\mathcal H_m$. If $H \in \mathcal H_m$, then $\mathcal H_m$ contains all $s$-hypergraphs $G \supset H$, such that $(V(H),G)$ is a cyclic $m$-extension; $\mathcal H_m$ also contains all $s$-hypergraphs $\tilde{H}$ such that $V(\tilde{H}) = V(H)$, $E(\tilde{H}) \supset E(H)$ and $\rho^{\max}(\tilde{H}) < \frac{m}{m(s - 1) - 1}$.

Note that for any $G \in \mathcal H_m$ there exists a sequence of hypergraphs 
$G_0 = (\{x\}, \varnothing) \subsetneq G_1 \ldots \subsetneq G_t \subseteq G$ such that $G_{i + 1}$ is a cyclic $m$-extension of $G_i$ for all $i \in \{0, \ldots, t - 1\}$. Let us call such sequence of hypergraphs \textit{$m$-decomposition} of $G$. 

Let $H \subset G$ be two subhypergraphs in a hypergraph $\Gamma$. The pair $(G, H)$ is \textit{cyclically $m$-maximal in $\Gamma$}, if there
are no cyclic $m$-extensions of $G$ in $\Gamma$ which are not cyclic $m$-extensions of $H$.
\\

\section{Proofs of theorems}
\indent

\textbf{Proof of Theorem 6}.
Let $m \in \mathbb N$, $\alpha = s - 1 - \frac{1}{2^{k-s-4}} + \frac{1}{2^{k - s - 4}m}$ and $p = n^{-\alpha}$.

Set 
\begin{align*}
& D_1(x_1, x_2) = (x_1 = x_2) \lor (\exists x_3 \ldots \exists x_s \; N(x_1, \ldots, x_s)), \\
& D_{i}(x_1, x_2) = \exists x_3 \mbox{ } (D_{\lfloor i/2 \rfloor}(x_1, x_3) \land D_{\lceil i/2 \rceil}(x_3, x_2)), \quad i > 1, \\
& D^{=}_1(x_1, x_2) = D_1(x_1, x_2) \land (x_1 \neq x_2), \\
& D^{=}_i(x_1, x_2) = D_i(x_1, x_2) \land (\neg (D_{i - 1}(x_1, x_2))), \quad i > 1.
\end{align*} 
The quantifier depths of $D_i$ and $D^{=}_i$ equal $\lceil \log_2 i \rceil + s - 2$. The formula $D_i(x_1, x_2)$ expresses the property that the distance between $x_1$ and $x_2$ is at most $i$, $D^{=}_i(x_1, x_2)$ expresses the property that the distance between $x_1$ and $x_2$ is exactly $i$. Moreover, set $D^{=}_{i, j}(x, y, z) = D^{=}_i(x, z) \land D^{=}_j(z, y)$
and $N_{i}(x, y) = \{z: D^{=}_{i, i}(x, y, z)\}$.

Let $L$ be a first-order property which is expressed by the formula $\exists a \exists b \mbox{ } Q(a, b)$ with quantifier depth $k$, where $Q(a, b) = Q_1(a, b) \land Q_2(a, b)$,
\begin{multline*}
Q_1(a, b) = D^{=}_{2^l}(a, b) \land (\neg (\exists u_1 \exists u_2 \mbox{ } ((u_1 \neq u_2) \land D^{=}_{2^{l-1}, 2^{l-1}}(a, b, u_1) \land D^{=}_{2^{l-1}, 2^{l-1}}(a, b, u_2) \land \\ 
\land (R_1(a, u_1, u_2) \lor R_2(a, u_1, u_2) \lor R_1(b, u_1, u_2) \lor R_2(b, u_1, u_2))))),
\end{multline*}
\begin{multline*}
Q_2(a, b) = \neg (\exists c \exists z_1 \exists z_2 \mbox{ } ((z_1 \neq z_2) \land (\neg D_{2^l}(a, z_1)) \land 
(\neg D_{2^l}(b, z_1)) \land (\neg D_{2^l}(a, z_2)) \land (\neg D_{2^l}(b, z_2)) \land  \\ \land
\forall u \mbox{ } ((D^{=}_{2^{l-1}, 2^{l-1}}(a, b, u) \land (u \neq c))  \Rightarrow (D^{=}_{2^l}(u, z_1) \lor D^{=}_{2^l}(u, z_2)))) ),
\end{multline*}

$$
R_1(a, u_1, u_2) = \exists x \mbox{ } \Bigl( \bigvee_{i=1}^{2^{l-1} - 1} (D^{=}_{i, 2^{l-1} - i}(u_1, a, x) \land D^{=}_i(u_2, x)) \Bigr),
$$
 
$$
R_2(a, u_1, u_2) = \exists x_1 \exists x_2 \mbox{ } (D^{=}_{2^{l-1} - 1}(u_1, x_1) \land D^{=}_{2^{l-1} - 1}(u_2, x_2) \land (\exists x_3 \ldots \exists x_{s - 1} \mbox{ } N(a, x_1, \ldots, x_{s - 1})))
$$

and $l = k - s - 4$. 

The predicate $Q_1(a, b)$ expresses the property that the distance between vertices $a$ and $b$ equals $2^l$ and there do not exist two edge-intersecting paths with length $2^{l-1}$ which connect the vertex $a$ (or, respectively, the vertex $b$) and two distinct vertices from the set $N_{2^{l-1}}(a, b)$.

Let $\tilde \Omega_n$ be the set of all hypergraphs $\mathcal G$ from $\Omega_n$ which satisfy the following properties.
\begin{itemize}
\item[1)] For any strictly balanced pair $(G, H)$ such that $\rho(G, H) < \frac{1}{\alpha}$ and $v(G) \le 2^{l + 2} s m$,  any $v(H)$-tuple has a strict $(G, H)$-extension in $\mathcal G$ which is $(K, T)$-maximal for all $(K, T) \in \mathcal K_{2^l s}$.
\item[2)] For any hypergraph $G$ with $\rho^{\max}(G) > \frac{1}{\alpha}$ and $v(G) \le 2^{l + 2} s m$, there is no copy of $G$ in $\mathcal G$. 
\end{itemize}

Theorems 9 and 11 imply that $\Pr[G^s(n, n^{-\alpha}) \in \tilde \Omega_n] \to 1$, $n \to \infty$.

Denote by $d_{\mathcal G}(x, y)$ the distance between vertices $x, y$ in a hypergraph $\mathcal G$.

Suppose that a hypergraph $\mathcal G \in \tilde \Omega_n$ satisfies $L$. Let $a, b$ be vertices such that $Q(a, b)$ is true.
Let $\chi = |N_{2^{l-1}}(a, b)|$ and $N_{2^{l-1}}(a, b) = \{x_1, \ldots, x_{\chi}\}$. 
Let $X$ be the union of $\chi$ paths with length $2^l$ and middle vertices $x_1, \ldots, x_{\chi}$ respectively which connect $a$ and $b$ in $\mathcal G$. 
Let us prove that $\chi \ge 2m$. Suppose that $\chi < 2m$. By property 1) from the definition of $\tilde \Omega_n$ in $\mathcal G$ there exist distinct vertices $z_1, z_2$ such that for any $i \in \{1, \ldots, \lfloor \chi/2 \rfloor \}$
the predicate $D^{=}_{2^{l}}(x_i, z_1)$ is true and for any $i \in \{\lfloor \chi/2 \rfloor + 1, \ldots, 2\lfloor \chi/2 \rfloor\}$ the predicate $D^{=}_{2^{l}}(x_i, z_2)$ is true, 
$d_{\mathcal G}(a, z_i) > 2^l$ and $d_{\mathcal G}(b, z_i) > 2^l$ for any $i \in \{1, 2\}$, 
and there exist $\lfloor \chi/2 \rfloor$ paths $P_1, \ldots, P_{\lfloor \chi/2 \rfloor}$ connecting $z_1$ and 
$x_1, \ldots, x_{\lfloor \chi/2 \rfloor}$ respectively such that for any distinct $i, j \in \{1, \ldots, \lfloor \chi/2 \rfloor \}$ equality $V(P_i) \cap V(P_j) = \{z_1\}$ holds and there exist $\lfloor \chi/2 \rfloor$ paths $Q_1, \ldots, Q_{\lfloor \chi/2 \rfloor}$ connecting $z_2$ and 
$x_{\lfloor \chi/2 \rfloor + 1}, \ldots, x_{2 \lfloor \chi/2 \rfloor}$ respectively such that for any distinct $i, j \in \{1, \ldots, \lfloor \chi/2 \rfloor \}$ equality $V(Q_i) \cap V(Q_j) = \{z_2\}$ holds. Indeed, in this case the pairs 
$(X \cup P_1 \cup \ldots \cup P_{\lfloor \chi/2 \rfloor}, X)$ and 
$(X \cup P_1 \cup \ldots \cup P_{\lfloor \chi/2 \rfloor} \cup Q_1 \cup \ldots \cup Q_{\lfloor \chi/2 \rfloor}, 
X \cup P_1 \cup \ldots \cup P_{\lfloor \chi/2 \rfloor})$ are strictly balanced and their densities equal
$$
\frac{2^l \lfloor \chi/2 \rfloor}{\lfloor \chi/2 \rfloor (2^l (s - 1) - 1) + 1} = \frac{1}{s - 1 - \frac{1}{2^l} + \frac{1}{\lfloor \chi/2 \rfloor 2^l}} < \frac{1}{s - 1 - \frac{1}{2^l} + \frac{1}{m 2^l}}= \frac{1}{\alpha}.
$$
This contradicts the truth of the predicate $Q_2(a, b)$. Therefore, $\chi \ge 2m$. Let us prove that $\chi \le 2m$. Suppose that $\chi > 2m$. 
Let $\tilde X$ be the union of $2m+1$ paths with length $2^l$ and middle vertices $x_1, \ldots, x_{2m+1}$ respectively which connect $a$ and $b$ in $\mathcal G$. 
Then the density of the subhypergraph $\tilde X$ is at least $$\frac{2^l (2m+1)}{(2m+1) (2^l (s - 1) - 1) + 2} > \frac{1}{\alpha}.$$ This contradicts property 2) from the definition of $\tilde \Omega_n$. Therefore, $\chi = 2m$.

Let $H$ be the union of $2m$ non-intersecting paths with length $2^l$ (such that consecutive edges of the paths intersect in one vertex) connecting vertices $a, b$. Let $x_1, \ldots, x_{2m}$ be the middle vertices of these paths. Let $G \supset H$ be an $s$-hypergraph obtained from $H$ by adding a vertex $z$ and $m$ non-intersecting paths with length $2^l$ connecting $z$ and vertices $x_1, \ldots, x_m$ respectively. Then $H$ is strictly balanced, the pair $(G, H)$ is strictly balanced and $\rho(H) = \rho(G, H) = \frac{1}{\alpha}$.
Let $\tilde L$ be the property that in $\mathcal G$ there exists a copy of $H$ such that no copy of $G$ contains it.

Let us show that if the hypergraph $\mathcal G \in \tilde \Omega_n$ satisfies $L$ then it satisfies the property $L_H$. Suppose that $\mathcal G$ satisfies $L$. Then there exist vertices $a, b$ satisfying $Q(a, b)$. As we have shown, 
$|N_{2^{l-1}}(a, b)| = 2m$. The property 2) from the definition of $\tilde \Omega_n$ implies that the hypergraph 
$X$ is isomorphic to $H$, so $\mathcal G$ contains a copy of $H$.

Let us prove that if the hypergraph $\mathcal G \in \tilde \Omega_n$ satisfies $\tilde L$ then it satisfies $L$.
Suppose that $\mathcal G$ satisfies $\tilde L$. Let $X_0 \subset \mathcal G$ be a copy of $H$ such that no copy of $G$ contains it. Let $a, b$ be the endpoints of the paths of $X_0$. Then $|N_{2^{l-1}}(a, b)| = 2m$, $Q_1(a, b)$ is true and
$d_{\mathcal G}(x, y) = 2^l$ for any distinct vertices $x, y \in N_{2^{l-1}}(a, b)$. 
Indeed, otherwise $\mathcal G$ contains a subhypergraph with at most $2^{l + 2} s m$ vertices and density greater than 
$1/\alpha$ which contradicts the definition of $\tilde \Omega_n$. Suppose that $Q(a, b)$ is false.
Since $Q_1(a, b)$ is true, there exist distinct vertices $z_1, z_2$ such that the predicate $D^{=}_{2^l}(\cdot, z_1) \lor D^{=}_{2^l}(\cdot, z_2)$ is true for all vertices from $N_{2^{l-1}}(a, b)$ except at most one and the predicates $D_{2^l}(a, z_i)$ and $D_{2^l}(b, z_i)$ are false for any $i \in \{1, 2\}$. Then without loss of generality we can assume that in $\mathcal G$ there exist paths $P_1, \ldots, P_m$ with length $2^l$ which connect vertex $z_1$ with vertices $x_1, \ldots, x_m$ respectively. Suppose that for some $i \in \{1, \ldots, m-1\}$ we have $P_{i+1} \subseteq P_1 \cup \ldots \cup P_i \cup X_0$. Let $P_{i + 1} \colon x_{i+1} \in {\bf e_1} \to \ldots \to {\bf e_{2^l}} \ni z_1$, where ${\bf e_1}, \ldots, {\bf e_{2^l}}$ are the edges of the path $P_{i+1}$ considered in the order from $x_{i+1}$ to $z_1$. Let $t$ be the minimal number such that 
${\bf e_t} \notin X_0$. Then ${\bf e_t} \in E(P_j)$ for some $j < i + 1$. Since $d_{\mathcal G}(x_{i+1}, z) = d_{\mathcal G}(x_j, z) = 2^l$, the edge ${\bf e_t}$ is also the $t$-th edge in the path $P_j$ considered in the order from $x_j$ to $z_1$. Therefore, $d_{\mathcal G}(x_{i+1}, x_j) \le 2t - 1$. Since $d_{\mathcal G}(x_{i+1}, x_j) = 2^l$, we obtain that $t \ge 2^{l-1} + 1$. 
This implies that one of the vertices $a, b$ belongs to the edge ${\bf e_{2^{l-1}}}$ and hence $\min(d_{\mathcal G}(a, z_1), d_{\mathcal G}(b, z_1)) \le 2^{l - 1} + 1$ which contradicts the falsity of the predicates $D_{2^l}(a, z_1)$ and $D_{2^l}(b, z_1)$. Therefore, $P_{i+1}  \not \subseteq P_1 \cup \ldots \cup P_i \cup X_0$ for all $i \in \{1, \ldots, m-1\}$. 
Consider the sequence of hypergraphs $X_0$, $X_1 = X_0 \cup P_1$, 
$X_2 = X_0 \cup P_1 \cup P_2$,
\ldots, $X_m = X_0 \cup P_1 \cup \ldots \cup P_m$. For any $i \in \{1, \ldots, m\}$ the hypergraph $X_i$ is obtained from the hypergraph $X_{i - 1}$ by adding $n_i$ vertices and
$e_i$ edges, where $e_i \le 2^l$ for any $i \in \{1, \ldots, m\}$ and $n_i \le e_i(s - 1) - 1$ for any $i \in \{2, \ldots, m\}$, $n_1 \le e_1(s - 1)$. Therefore, 
\begin{multline*}
1/\rho(X_m) \le \frac{2m (2^l (s - 1) - 1) + 2 + (e_1 + \ldots + e_m)(s - 1) - m + 1}{2m 2^l + e_1 + \ldots + e_m} = \\ 
= s - 1 - \frac{3(m-1)}{2m 2^l + e_1 + \ldots + e_m} \le \alpha.
\end{multline*}
Equalities hold if and only if $e_i = 2^l$ for all $i \in \{1, \ldots, m\}$ and $n_i = 2^l (s - 1) - 1$ for all $i \in \{2, \ldots, m\}$,
$n_1 = 2^l (s - 1)$. 
Therefore, by the definition of $\tilde \Omega_n$ these equalities hold and $X_m$ is isomorphic to $G$. This contradicts the property $\tilde L$. Therefore, $Q(a, b)$ is true, and hence $\mathcal G$ satisfies $L$.

By Proposition 1, there exists $\lim_{n \to \infty} \Pr[G^s(n, n^{-\alpha})\models \tilde L] = c_1 \in (0, 1)$. 
By Theorem 10, there exists $\lim_{n \to \infty} \Pr[G^s(n, n^{-\alpha})\models L_H] = c_2 \in (0, 1)$.
Since $\lim_{n \to \infty} \Pr[G^s(n, n^{-\alpha}) \in \tilde \Omega_{n}] = 1$, it follows from the above that 
$\lim \inf_{n \to \infty} \Pr[G^s(n, n^{-\alpha}) \models L] \ge  c_1$,  
$\lim \sup_{n \to \infty} \Pr[G^s(n, n^{-\alpha}) \models L] \le  c_2$. 
Letting $m \to \infty$ we obtain that $s - 1 - \frac{1}{2^{k - s - 4}} \in (S_k)'$.
\\

\textbf{Proof of Theorem 7}. 

For an $s$-hypergraph $G$, let $d_G(x, \tilde G) = \min_{y \in \tilde G} d_G(x, y)$ denote the distance between a vertex $x \in V(G)$ and a subhypergraph $\tilde G \subset G$.
Denote by $d_G(\tilde G_1, \tilde G_2) = \min_{x \in V(G_1)} d_G(x, G_2)$ the distance between $\tilde G_1, \tilde G_2 \subset G$.

By Proposition 1 from \cite{zero_one} there exists $\eta(1/\alpha) \in \mathbb N$ such that 
any hypergraph $H_1 \in \mathcal H_{2^{k-s+1}}$ with $v(H_1) \ge \eta(1/\alpha)$ contains a subhypergraph $H_2 \in \mathcal H_{2^{k-s+1}}$ with $v(H_2) \le \eta(1/\alpha)$ and density $\rho(H_2) > 1/\alpha$. 
Set $n(\alpha) = \eta(1/\alpha) + (k - s + 1) s 2^{k - s + 1}$. 

Let $\tilde \Omega_n$ be the set of all hypergraphs $G$ from $\Omega_n$ which satisfy the following properties.
\begin{itemize}
\item[1)] Let $H$ be an $s$-hypergraph with at most $n(\alpha)$ vertices. 
If $\rho^{\max}(H) > 1/\alpha$, then $G$ does not contain a subhypergraph isomorphic to $H$. 
If $\rho^{\max}(H) < 1/\alpha$, then $G$ contains an induced subhypergraph which is isomorphic to $H$ and $(K_2, K_1)$-maximal in $G$ for all $(K_2, K_1) \in \mathcal K_{2^{k - s + 1} s}$.  
\item[2)] For any $\alpha$-safe pair $(H_2, H_1)$ with $v(H_2) \le n(\alpha)$ and for any $G_1 \subset G$ with $v(G_1) = v(H_1)$, there exists a subhypergraph $G_2 \subset G$ such that $G_2$ is a strict $(H_2, H_1)$-extension of $G_1$
and $(G_2, G_1)$ is $(K_2, K_1)$-maximal in $G$ for all $(K_2, K_1) \in \mathcal K_{2^{k - s + 1} s}$. 
\end{itemize}

Theorems 9, 11 imply that $\Pr[G^s(n, n^{-\alpha}) \in \tilde \Omega_n] \to 1$ as $n \to \infty$. The statement of the theorem follows from the existence of a winning strategy for Duplicator in $\EHR(G, H, k)$ for all pairs $(G, H)$ such that $G, H \in \tilde \Omega_n$.

Let $G, H \in \tilde \Omega_n$. Let $X_r$ and $Y_r$ be the hypergraphs chosen in the $r$-th round by Spoiler and Duplicator respectively. We denote vertices which are chosen in the first $r$ rounds in $X_r$ and $Y_r$ by $x^1_r, \ldots, x^r_r$ and $y^1_r, \ldots, y^r_r$.
Let us describe Duplicator's strategy by induction.

Let $r$ rounds be finished, where $1 \le r \le k - s + 2$. Let $l \in \{1, \ldots, r\}$. Let $\tilde X^1_r, \ldots, \tilde X^l_r \subset X_r$ and $\tilde Y^1_r, \ldots, \tilde Y^l_r \subset Y_r$ be subhypergraphs of $X_r$ and $Y_r$ respectively. We say that $\tilde X^1_r, \ldots, \tilde X^l_r$ and $\tilde Y^1_r, \ldots, \tilde Y^l_r$ are \textit{$(k, r, l)$-regular equivalent} in $(X_r, Y_r)$, if the following properties hold.
\begin{itemize}
\item[(I)] $x^1_r, \ldots, x^r_r \in V(\tilde X^1_r \cup \ldots \cup \tilde X^l_r)$, $y^1_r, \ldots, y^r_r \in V(\tilde Y^1_r \cup \ldots \cup \tilde Y^l_r).$
\item[(II)] For any distinct $j_1, j_2 \in \{1, \ldots, l\}$, the inequalities $d_{X_r}(\tilde X^{j_1}_r, \tilde X^{j_2}_r) > 2^{k - r - s + 2}$,
$d_{Y_r}(\tilde Y^{j_1}_r, \tilde Y^{j_2}_r) > 2^{k - r - s + 2}$ hold.
\item[(III)] For any $j \in \{1, \ldots, l\}$, there is no cyclic $2^{k - r - s + 2}$-extension of $\tilde X^j_r$ in the hypergraph $X_r$
and there is no cyclic $2^{k - r - s + 2}$-extension of $\tilde Y^j_r$ in the hypergraph $Y_r$.
\item[(IV)] Cardinalities of the sets $V(\tilde X^1_r \cup \ldots \cup \tilde X^l_r)$ and $V(\tilde Y^1_i \cup \ldots \cup \tilde Y^l_r)$ are at most $\eta(\rho) + (r - 1) s 2^{k - s + 1}$.
\item[(V)] The hypergraphs $\tilde X^j_r$ and $\tilde Y^j_r$ are isomorphic for any $j \in \{1, \ldots, l\}$ and there exists a corresponding isomorphism (one for all these pairs of hypergraphs) which maps the vertices $x^i_r$ to the vertices $y^i_r$ for all $i \in \{1, \ldots, r\}$.
\end{itemize}

Two hypergraphs $\tilde X^1_r$ and $\tilde Y^1_r$ are called $(k, r)$-equivalent in $(X_r, Y_r)$ if
\begin{itemize}
\item[a)] properties (I), (IV) and (V) hold for $l = 1$,
\item[b)] there is no cyclic $(2^{k - r - s + 2} - 1)$-extension of $\tilde X^1_r$ (resp. $\tilde Y^1_r$)  in $X_r$ (resp. $Y_r$),
\item[c)] there is no second type cyclic $2^{k - r - s + 2}$-extension of $X_r|_{\{x^1_r, \ldots, x^r_r\}}$ (resp. $Y_r|_{\{y^1_r, \ldots, y^r_r\}}$) in $X_r \setminus (\tilde X^1_r \setminus X_r|_{\{x^1_r, \ldots, x^r_r\}})$ 
(resp. $Y_r \setminus (\tilde Y^1_r \setminus Y_r|_{\{y^1_r, \ldots, y^r_r\}})$),
\end{itemize}

The main idea of Duplicator's strategy is the following. Duplicator should play in such way
that for some $r \in \{1, \ldots, k - s + 1\}$ and $l \in \{1, \ldots, r\}$ in the hypergraphs $X_r$, $Y_r$ $(k, r, l)$-regular equivalent collections of subhypergraphs are constructed. In the first round, Duplicator must use the
strategy $S_1$ which is described in the next section. After the $r$-th round, $r \in \{1, \ldots, k - s\}$, if
$(k, r, l)$-regular equivalent collections are not constructed, then Duplicator can find $(k, r)$-equivalent hypergraphs and in the $(r + 1)$-th round he must use the strategy $S_{r+1}$. Strategy SF is described in \cite{zero_one} and is used by Duplicator in the $(r + 1)$-th round, $r \ge 1$, if and only if after the $r$-th round for some $l \in \{1, \ldots, r\}$ $(k, r, l)$-regular equivalent collections of hypergraphs in 
$(X_r, Y_r)$ are constructed. In \cite{zero_one} it is proved that Duplicator wins, when he uses the strategy SF.\\

\textbf{Strategy $S_1$}

Consider the first round and two possibilities to choose the first vertex by Spoiler. 

Let in $X_1$ there is no cyclic $2^{k - s + 1}$-extension of $(\{x^1_1\}, \varnothing)$. Then Duplicator chooses a vertex $y^1_1 \in V(Y_1)$ such that there are no cyclic $2^{k - s + 1}$-extensions of $(\{y^1_1\}, \varnothing)$ in $Y_1$ (such a vertex exists because $Y_1 \in \tilde \Omega_n$). Set $\tilde X^1_1 = (\{x^1_1\}, \varnothing)$, $\tilde Y^1_1 = (\{y^1_1\}, \varnothing)$. Then $\tilde X^1_1$ and $\tilde Y^1_1$ are $(k, 1, 1)$-regular equivalent in $(X_1, Y_1)$. In this case, in the second round Duplicator exploits the strategy SF.

Let in $X_1$ there exists at least one cyclic $2^{k - s + 1}$-extension of $(\{x^1_1\}, \varnothing)$. Let us prove that there exists a sequence of hypergraphs $G_1, G_2, \ldots, G_t$ such that

\begin{itemize}
\item[a)] for any $i \in \{1, \ldots, s - 1\}$ the hypergraph $G_{i + 1}$ is a cyclic $2^{k - s + 1}$-extension of the hypergraph $G_i$ in $X_1$, $G_1$ is a cyclic $2^{k - s + 1}$-extension of $(\{x^1_1\}, \varnothing)$,

\item[b)] $G_t$ is an induced subhypergraph of $X_1$, 

\item[c)] there are no cyclic $2^{k - s + 1}$-extensions of $G_t$ in $X_1$,

\item[d)] either $\rho(G_t) < 1/\alpha$ or there exists $i \in \{1, \ldots, t - 1\}$ such that the hypergraph $G_{i + 1}$ is a cyclic $2^{k - s + 1}$-extension of the hypergraph $G_i$ and there are no cyclic $(2^{k - s + 1} - 1)$-extensions of $G_i$ in $X_1$. 

\end{itemize}

Let us prove the existence of such a sequence. Obviously, there exists a sequence $G_0 \subset  G_1 \subset \ldots \subset G_i$ with the following properties. First, $G_0 = \{x^1_1\}, \varnothing)$ and $G_j$ is a cyclic $(2^{k - s + 1}-1)$-extension of the hypergraph $G_{j - 1}$ for any $j \in \{1, \ldots, i\}$. Second, $j = i$ is the first number such that $G_j$ has no cyclic $(2^{k - s + 1}-1)$-extensions in $X_1$. Let us add cyclic $2^{k - s + 1}$-extensions to the hypergraph $G_i$ (each next hypergraph is a cyclic $2^{k-s+1}$-extension of the previous one) until there are no cyclic $2^{k - s + 1}$-extensions of the final hypergraph in $X_1$. 
We get the sequence of hypergraphs $G_1, \ldots, G_t$, which follows Properties a) and c) (in addition, the inequality 
$t \le 2^{k - s + 1} b + 1$ holds, because the density of $G_t$ is greater than $1/\alpha$, if $t = 2^{k - s + 1} b + 2$, this contradicts Property 1)).

Let us prove that $G_t$ is an induced subhypergraph of $X_1$. 
Set $e_i = e(G_i, G_{i - 1})$, $v_i = v(G_i, G_{i - 1})$ for all $i \in \{1, \ldots, t\}$. Set $e_0 = e(X_1|_{V(G_t)}) - e(G_t)$. 
Note that $e_i \le 2^{k - s + 1}$, $v_i \le e_i (s - 1) - 1$. Then 
$$
\frac{1}{\rho(X_1|_{V(G_t)})} \le \frac{(e_1 + \ldots + e_t) (s - 1) - t + 1}{e_1 + \ldots + e_t + e} = s - 1 -
\frac{(s - 1)e + t - 1}{e_1 + \ldots + e_t + e} = s - 1 - \frac{1}{\tau},
$$
where $\tau = \frac{e_1 + \ldots + e_t + e}{(s - 1)e + t - 1}$. We have
$$
\tau - 2^{k - s + 1} = \frac{(e_1 - 2^{k - s +1}) + \ldots + (e_t - 2^{k - s + 1}) + e - 2^{k - s + 1} ((s - 1)e - 1)}{(s - 1)e + t - 1}
$$
which is less than $0$ if $e \ge 1$ and $s \ge 3$. Thus $\frac{1}{\rho(X_1|_{V(G_t)})} < \alpha$ if $e \ge 1$. Therefore, $e = 0$ and $G_t$ is an induced subhypergraph of $X_1$.

Let us show that $v_i = e_i (s - 1) - 1$ for all $i \in \{1, \ldots, t\}$. Suppose that there exists $i \in \{1, \ldots, t\}$ such that $v_i \le e_i (s - 1) - 2$. Then
$$
\frac{1}{\rho(G_t)} \le \frac{(e_1 + \ldots + e_t) (s - 1) - t}{e_1 + \ldots + e_t} = s - 1 - \frac{1}{(e_1 + \ldots + e_t)/t} 
\le s - 1 - \frac{1}{2^{k - s + 1}} < \alpha. 
$$
Therefore, $v_i = e_i (s - 1) - 1$ for all $i \in \{1, \ldots, t\}$.

Let us prove that the hypergraph $G_t$ is strictly balanced. Let $\tilde G$ be an arbitrary proper subhypergraph in $G_t$. Denote $\tilde G_1 = G_1 \cap \tilde G$. If $\tilde G_1 \neq G_1$, then $e(\tilde G \cup G_1, \tilde G) \le 2^{k - s + 1}$,
$v(\tilde G \cup G_1, \tilde G) \le e(\tilde G \cup G_1, \tilde G) (s - 1) - 1$. Therefore, the density of the hypergraph $\tilde G \cup G_1$ is at least
$$
\frac{e(\tilde G) + e(\tilde G \cup G_1, \tilde G)}{v(\tilde G) + e(\tilde G \cup G_1, \tilde G) (s - 1) - 1} > \min \Bigl(\frac{e(\tilde G)}{v(\tilde G)}, \frac{1}{s - 1 - \frac{1}{e(\tilde G \cup G_1, \tilde G)}}\Bigr) = \rho(\tilde G).
$$
In the same way, it can be proved that $\rho(G_t) \ge \rho(\tilde G \cup G_{t - 1}) \ge \ldots \ge \rho(\tilde G \cup G_1) \ge \rho(\tilde G)$, where at least one of the inequalities is strict, because $\tilde G$ is a proper subhypergraph in $G_t$. Therefore, the hypergraph $G_t$ is strictly balanced.

If $\rho(G_t) < \frac{1}{\alpha}$, then set $\tilde X^1_1 = G_t$. By the definition of the set $\tilde \Omega_n$ 
the hypergraph $Y_1$ has a subhypergraph $\tilde Y^1_1$ which is isomorphic to $\tilde X^1_1$ and $(K, T)$-maximal for any pair $(K, T)$ such that $v(K) \le 2^{k - s + 1} s$ and $f_{\alpha}(K, T) < 0$. Let $\varphi \colon \tilde X^1_1 \to \tilde Y^1_1$ be an isomorphism. Then Duplicator chooses the vertex $y^1_1 = \varphi(x^1_1)$. By the construction of the hypergraphs $\tilde X^1_1$ and $\tilde Y^1_1$, they
do not have cyclic $2^{k - s + 1}$-extensions in $X_1$ and $Y_1$ respectively. Therefore, the hypergraphs $\tilde X^1_1$ and $\tilde Y^1_1$ are $(k, 1, 1)$-regular equivalent in $(X_1, Y_1)$. In the second round Duplicator exploits the strategy SF.

Let $\rho(G_t) = \frac{1}{\alpha}$. Then
$$
1/\rho(G_t) = s - 1 - \frac{1}{(e_1 + \ldots + e_t)/(t - 1)} = s - 1 - \frac{1}{2^{k - s + 1} + a/b}.
$$
Since $a/b$ is an irreducible fraction, $t \ge b + 1$. Let us show that there exists $i \in \{0, \ldots, t - 1\}$ such that $G_{i + 1}$ is not a cyclic $(2^{k - s + 1} - 1)$-extension of $G_i$. Indeed, otherwise
\begin{multline*}
1/\rho(G_t) \le s - 1 - \frac{1}{2^{k - s + 1} - 1 + \frac{2^{k - s + 1} - 1}{t - 1}} \le
s - 1 - \frac{1}{2^{k - s + 1} - 1 + \frac{2^{k - s + 1} - 1}{b}} = \\
= s - 1 - \frac{1}{2^{k - s + 1} + \frac{2^{k - s + 1} - b - 1}{b}} < \alpha,
\end{multline*}
since $a \ge 2^{k - s + 1} - b$. 
Since $G_t$ is strictly balanced, $\rho^{\max}(G_i) < 1/\alpha$. As $Y_1 \in \tilde \Omega_n$, in $Y_1$ there exists a subhypergraph $\tilde Y^1_1$ which is isomorphic to $\tilde X^1_1 := G_i$ and $(K, T)$-maximal for any pair $(K, T)$ such that $v(K) \le 2^{k - s + 1} s$ and $f_{\alpha}(K, T) < 0$.
Let $\varphi \colon \tilde X^1_1 \to \tilde Y^1_1$ be an isomorphism. Then Duplicator chooses the vertex $y^1_1 = \varphi(x^1_1)$. By the construction of the hypergraphs $\tilde X^1_1$ and $\tilde Y^1_1$, they
are $(k, 1)$-equivalent in $(X_1, Y_1)$. Therefore, in the second round Duplicator exploits the strategy $S_2$.\\

\textbf{Strategy $S_{r+1}$}

Let after the $r$-th round, $r \in \{1, \ldots, k - s + 1\}$, there exist hypergraphs $\tilde X^1_r, \tilde Y^1_r$ which are $(k, r)$-equivalent in
$(X_r, Y_r)$. Let $\varphi \colon \tilde X^1_r \to \tilde Y^1_r$ be an isomorphism. In the $(r+1)$-th round, Spoiler chooses a vertex $x^{r+1}_{r+1}$. If $X_{r+1} = X_r$, then set $\tilde X^1_{r + 1} = \tilde X^1_r$, $\tilde Y^1_{r + 1} = \tilde Y^1_r$.
Otherwise, set  $\tilde X^1_{r + 1} = \tilde Y^1_r$, $\tilde Y^1_{r + 1} = \tilde X^1_r$.

Let $x^{r+1}_{r+1} \in V (\tilde X^1_{r + 1})$. Duplicator chooses the vertex $y^{r+1}_{r+1} = \varphi(x^{r+1}_{r+1})$,
if $X_{r+1} = X_r$, and the vertex $y^{r+1}_{r+1} = \varphi^{-1}(x^{r+1}_{r+1})$,
if $X_{r+1} = Y_r$. As in $X_r, Y_r$ there are no cyclic $2^{k - r - s + 1}$-extensions of the
hypergraphs $\tilde X^1_r, \tilde Y^1_r$ respectively (by the definition of $(k, r)$-equivalence), the hypergraphs 
$\tilde X^1_{r + 1}, \tilde Y^1_{r + 1}$ are $(k, r + 1, 1)$-regular equivalent in $(X_{r+1}, Y_{r+1})$. Therefore, in the $(r + 2)$-th round Duplicator exploits the strategy SF.

Let $x^{r+1}_{r+1} \notin V (\tilde X^1_{r + 1})$. Consider several cases. 

Suppose first that $r < k - s + 1$.

Let $d_{X_{r+1}}(\tilde X^1_{r+1}, x^{r+1}_{r+1}) > 2^{k - r - s + 1}$. If in $X_{r + 1}$ there are no cyclic $2^{k - r - s + 1}$-extensions of the hypergraph $(\{x^{r+1}_{r+1}\}, \varnothing)$, then set $\tilde X^2_{r+1} = (\{x^{r+1}_{r+1}\}, \varnothing)$. By Property 2) of the hypergraph $Y_{r+1}$, it has a vertex $y^{r+1}_{r+1}$ such that $d_{Y_{r+1}}(\tilde Y^1_{r+1}, y^{r+1}_{r+1}) = 2^{k - r - s + 1} + 1$ and there are no cyclic $2^{k - r - s + 1}$-extensions of $(\{y^{r+1}_{r+1}\}, \varnothing)$ in $Y_{r + 1}$.
Set $\tilde Y^2_{r+1} = (\{y^{r+1}_{r+1}\}, \varnothing)$. If there is exactly one cyclic $2^{k - r - s + 1}$-extension of $(\{x^{r+1}_{r+1}\}, \varnothing)$, then we denote it by $\tilde X^2_{r + 1}$ (the hypergraph $(\{x^{r+1}_{r+1}\}, \varnothing)$ has at most one cyclic $2^{k - r - s + 1}$-extension in $X_{r+1}$). Let $d_{X_{r+1}}(\tilde X^1_{r+1}, \tilde X^2_{r+1}) > 2^{k - r - s + 1}$. By Property 2) of the hypergraph $Y_{r+1}$, it has
a vertex $y^{r+1}_{r+1}$ and a subhypergraph $\tilde Y^2_{r+1}$ such that $d_{Y_{r+1}}(\tilde Y^1_{r+1}, \tilde Y^2_{r+1}) = 2^{k - r - s + 1} + 1$, pairs $(\tilde Y^2_{r+1}, (\{y^{r+1}_{r+1}\}, \varnothing))$ and $(\tilde X^2_{r+1}, (\{x^{r+1}_{r+1}\}, \varnothing))$ are isomorphic, and there are no cyclic $2^{k - r - s + 1}$-extensions of $\tilde Y^2_{r+1}$ in $Y_{r+1}$.
The property of $(k, r)$-equivalence of the hypergraphs $\tilde X^1_r$, $\tilde Y^1_r$ in $(X_r, Y_r)$ implies non-existence of cyclic
$2^{k - r - s + 1}$-extensions of $\tilde X^1_r$ and $\tilde Y^1_r$ in $X_r$ and $Y_r$ respectively. Obviously, in all the considered cases
the collections $\tilde X^1_{r + 1}$, $\tilde X^2_{r + 1}$ and $\tilde Y^1_{r + 1}$, $\tilde Y^2_{r + 1}$ are $(k, r + 1, 2)$-regular equivalent in $(X_{r+1}, Y_{r+1})$. Thus, in the $(r + 2)$-th round Duplicator exploits the strategy SF.
Let $d_{X_{r+1}}(\tilde X^1_{r+1}, \tilde X^2_{r+1}) \le 2^{k - r - s + 1}$. The property of $(k, r)$-equivalence of the hypergraphs $\tilde X^1_r$, $\tilde Y^1_r$ in $(X_r, Y_r)$ implies $d_{X_{r+1}}(\tilde X^1_{r+1}, \tilde X^2_{r+1}) = 2^{k - r - s + 1}$
and non-existence of cyclic $(2^{k - r - s + 1}-1)$-extensions of $(\{x^{r+1}_{r+1}\}, \varnothing)$ in $X_{r+1}$. 
By Property 2) of the hypergraph $Y_{r+1}$, it has a vertex $y^{r+1}_{r+1}$ such that $d_{Y_{r+1}}(\tilde Y^1_{r+1}, y^{r+1}_{r+1}) = 2^{k - r - s + 1} + 1$ and there are no cyclic $2^{k - r - s + 1}$-extensions of $(\{y^{r+1}_{r+1}\}, \varnothing)$ in $Y_{r + 1}$.
Set $\tilde X^1_{r + 1} :=\tilde X^1_{r + 1} \cup \{\{x^{r+1}_{r+1}\}, \varnothing\}$, $\tilde Y^1_{r + 1} :=\tilde Y^1_{r + 1} \cup \{\{y^{r+1}_{r+1}\}, \varnothing\}$. The hypergraphs $\tilde X^1_{r + 1}$ and $\tilde Y^1_{r + 1}$ are $(k, r+1)$-equivalent in 
$(X_{r+1}, Y_{r+1})$ and in the next round Duplicator exploits the strategy $S_{r+2}$.

Let $d_{X_{r+1}}(\tilde X^1_{r+1}, x^{r+1}_{r+1}) \le 2^{k - r - s + 1}$. Consider a minimal chain $L_X$ in $X_{r+1}$ which connects $x^{r+1}_{r+1}$ and $\tilde X^1_{r+1}$. Moreover, let such chain connect $x^{r+1}_{r+1}$ and one of the vertices $x^1_{r+1}, \ldots, x^r_{r+1}$ if such a minimal chain exists. 
By Property 2) of the hypergraph $Y_{r+1}$, there exists a vertex $y^{r+1}_{r+1}$ such that $d_{Y_{r+1}}(\tilde Y^1_{r + 1}, y^{r+1}_{r+1}) = d_{X_{r+1}}(\tilde X^1_{r + 1}, x^{r+1}_{r+1})$, there exists an isomorphism 
$L_X \cup \tilde X^1_{r+1} \to L_Y \cup \tilde Y^1_{r+1}$ which maps 
the vertices $x^1_{r+1}, \ldots, x^{r+1}_{r+1}$ to the vertices $y^1_{r+1}, \ldots, y^{r+1}_{r+1}$ respectively, and the pair 
$(L_Y \cup \tilde Y^1_{r+1}, \tilde Y^1_{r+1})$ is cyclically $2^{k - r - s + 1}$-maximal in $Y_{r+1}$, where $L_Y$ is a minimal chain which connects $y^{r+1}_{r+1}$ and $\tilde Y^1_{r+1}$ in $Y_{r+1}$. Obviously,
there are no cyclic $2^{k - r - s + 1}$-extensions of the hypergraph $\tilde Y^1_{r+1}$ in $Y_{r+1}$. If there are no cyclic $2^{k - r - s + 1}$-extensions of $L_X \cup \tilde X^1_{r + 1}$ in $X_{r+1}$, then set 
$\tilde X^1_{r+1} := L_X \cup \tilde X^1_{r+1}$ and $\tilde Y^1_{r+1} := \tilde Y^1_{r+1} \cup L_Y$. Obviously, the
hypergraphs $\tilde X^1_{r+1}$ and $\tilde Y^1_{r+1}$ are $(k, r+1, 1)$-regular equivalent in $(X_{r+1}, Y_{r+1})$. Therefore, in the next round Duplicator exploits the strategy SF. If there is a cyclic $2^{k - r - s + 1}$-extension of $L_X \cup \tilde X^1_{r + 1}$ in $X_{r+1}$, then $d_{X_{r+1}}(x^{r+1}_{r+1}, \tilde X^1_{r+1}) = 2^{k - r - s + 1}$ and there are no cyclic $(2^{k - r - s + 1} - 1)$-extensions of $L_X \cup \tilde X^1_{r + 1}$ in $X_{r+1}$. 
The property of $(k, r)$-equivalence of the hypergraphs $\tilde X^1_{r+1}, \tilde Y^1_{r+1}$ in 
$(X_{r+1}, Y_{r+1})$ implies that either
there are no second type cyclic $2^{k - r - s + 1}$-extension of $X_{r + 1}|_{\{x^1_{r+1}, \ldots, x^{r+1}_{r+1}\}}$ in $X_{r + 1} \setminus (\tilde X^1_{r+1} \setminus X_{r+1}|_{\{x^1_{r+1}, \ldots, x^{r+1}_{r+1}\}})$ or 
the chain $L_X$ could be chosen from a set of chains $\{L_X^j\}$ which connect vertex $x^{r+1}_{r+1}$ with some vertex $x^i_{r+1}$, where $i \in \{1, \ldots, r\}$.  
In the first case the hypergraphs $\tilde X^1_{r+1} := L_X \cup \tilde X^1_{r+1}$ and $\tilde Y^1_{r+1} := \tilde Y^1_{r+1} \cup L_Y$ are $(k, r+1)$-equivalent in 
$(X_{r+1}, Y_{r+1})$ and in the next round Duplicator exploits the strategy $S_{r+2}$.
Consider the second case. If in the $(r+2)$-th round Spoiler chooses a vertex from the chain $L_X^j$ for some $j$, then set $\tilde X^1_{r+1} := L_X^j \cup \tilde X^1_{r+1}$. Otherwise, set 
$\tilde X^1_{r+1} := L_X^1 \cup \tilde X^1_{r+1}$. Set $\tilde Y^1_{r+1} := \tilde Y^1_{r+1} \cup L_Y$. Then the hypergraphs $\tilde X^1_{r+1}$ and $\tilde Y^1_{r+1}$ satisfy properties a), b) from the definition of $(k, r+1)$-equivalence and the vertex $x^{r+2}_{r+2}$ does not belong to any second type cyclic $2^{k - r - s + 1}$-extension of $X_{r + 1}|_{\{x^1_{r+1}, \ldots, x^{r+1}_{r+1}\}}$ in $X_{r + 1} \setminus (\tilde X^1_{r+1} \setminus X_{r+1}|_{\{x^1_{r+1}, \ldots, x^{r+1}_{r+1}\}})$
which enables Duplicator to exploit the strategy $S_{r+2}$ in the $(r+2)$-th round.
 
Suppose that $r = k - s + 1$. If $d_{X_{r+1}}(\tilde X^1_{r+1}, x^{r+1}_{r+1}) > 1$, then set 
$\tilde X^2_{r+1} = (\{x^{r+1}_{r+1}\}, \varnothing)$. By Property 2) of the hypergraph $Y_{r+1}$, it has a vertex $y^{r+1}_{r+1}$ such that $d_{Y_{r+1}}(\tilde Y^1_{r+1}, y^{r+1}_{r+1}) = 2$. 
Set $\tilde Y^2_{r+1} = (\{y^{r+1}_{r+1}\}, \varnothing)$.
The property of $(k, r)$-equivalence of the hypergraphs $\tilde X^1_r$, $\tilde Y^1_r$ in $(X_r, Y_r)$ implies non-existence of cyclic $1$-extensions of $\tilde X^1_r$ and $\tilde Y^1_r$ in $X_r$ and $Y_r$ respectively. Then 
the collections $\tilde X^1_{r + 1}$, $\tilde X^2_{r + 1}$ and $\tilde Y^1_{r + 1}$, $\tilde Y^2_{r + 1}$ are $(k, r + 1, 2)$-regular equivalent in $(X_{r+1}, Y_{r+1})$. Thus, in the $(r + 2)$-th round Duplicator exploits the strategy SF.

If $d_{X_{r+1}}(\tilde X^1_{r+1}, x^{r+1}_{r+1}) = 1$ then consider an edge $e_X$ which connects $x^{r+1}_{r+1}$ and $\tilde X^1_{r + 1}$. Moreover, let such edge connect $x^{r+1}_{r+1}$ and one of the vertices $x^1_{r+1}, \ldots, x^r_{r+1}$ if such an edge exists. 
By Property 2) of the hypergraph $Y_{r+1}$, there exists a vertex $y^{r+1}_{r+1}$ such that $d_{Y_{r+1}}(\tilde Y^1_{r + 1}, y^{r+1}_{r+1}) =1$, there exists an isomorphism 
$e_X \cup \tilde X^1_{r+1} \to e_Y \cup \tilde Y^1_{r+1}$ which maps 
the vertices $x^1_{r+1}, \ldots, x^{r+1}_{r+1}$ to the vertices $y^1_{r+1}, \ldots, y^{r+1}_{r+1}$ respectively, and the pair 
$(e_Y \cup \tilde Y^1_{r+1}, \tilde Y^1_{r+1})$ is cyclically $1$-maximal in $Y_{r+1}$, where $e_Y$ is an edge which connects $y^{r+1}_{r+1}$ and $\tilde Y^1_{r+1}$ in $Y_{r+1}$. 
If there are no cyclic $1$-extensions of $e_X \cup \tilde X^1_{r + 1}$ in $X_{r+1}$, then the hypergraphs $e_X \cup \tilde X^1_{r + 1}$ and $e_Y \cup \tilde Y^1_{r + 1}$ are $(k, r + 1, 2)$-regular equivalent in $(X_{r+1}, Y_{r+1})$. Thus, in the $(r + 2)$-th round Duplicator exploits the strategy SF. If there is an edge $\tilde e_X$ which forms a cyclic $1$-extension of $e_X \cup \tilde X^1_{r + 1}$ in $X_{r+1}$ then the property of $(k, r)$-equivalence of the hypergraphs $\tilde X^1_{r+1}$, $\tilde Y^1_{r+1}$ implies that either $\tilde e_X$ contains at most one vertex from the set $\{x^1_{r+1}, \ldots, x^{r+1}_{r+1}\}$ 
or there exists a vertex $x^i_{r+1}$, where $i \in \{1, \dots, r\}$, such that $x^i_{r+1}, x^{r+1}_{r+1} \in e_X \cap \tilde e_X$. In the first case set $\tilde X^1_{r+1} := e_X \cup \tilde X^1_{r + 1}$ and $\tilde Y^1_{r+1} := e_Y \cup \tilde Y^1_{r + 1}$. Consider the second case. If in the $(r+2)$-th round Spoiler chooses a vertex from the edge $\tilde e_X$, then set $\tilde X^1_{r+1} := \tilde e_X \cup \tilde X^1_{r+1}$. Otherwise, set $\tilde X^1_{r+1} := e_X \cup \tilde X^1_{r + 1}$. 
Set $\tilde Y^1_{r+1} := e_Y \cup \tilde Y^1_{r + 1}$.
Then Duplicator wins using the strategy SF for the hypergraphs $\tilde X^1_{r + 1}$ and $\tilde Y^1_{r + 1}$ in the rounds 
$k - s + 3, \ldots, k$, since there do not exist vertices $u_1, \ldots, u_l \notin V(\tilde X^1_{r+1})$, where $l \in \{1, \ldots, s-2\}$, and vertices $v_1, \ldots, v_{s - l} \in V(\tilde X^1_{r+1})$ such that $\{u_1, \ldots, u_l, v_1, \ldots, v_{s - l}\} \in E(X_{r+1})$ and $|\{v_1, \ldots, v_{s-l}\} \cap \{x^1_{r+1}, \ldots, x^{r+1}_{r+1}\}| \ge 2$.
\\

\textbf{Proof of Theorem 8}. 
Let $a \in \mathbb N$, $a \le 2^{k - s + 1} - 3$ and $\alpha = s - 1 - \frac{1}{2^{k - s + 1} + a}$. 

Condiser first the case where $k \ge s + 2$. 
Let $a + 3 = a_1 + a_2$, where $a_1, a_2 \in \{1, \ldots, 2^{k-s}\}$.

Let $L$ be a first-order property which is expressed by the formula $\exists x_1 \mbox{ } (Q_1(x_1) \land Q_2(x_1))$ with quantifier depth at most $k$, where 
$$
Q_1(x_1) = \exists x_2 \mbox{ } (D^{=}_{a_1}(x_1, x_2) \land (\exists x_3 \mbox{ } (D^{=}_{2^{k-s-1}}(x_2, x_3) 
\land \ldots \land (\exists x_{k-s} \mbox{ } (D^{=}_{2^2}(x_{k-s-1}, x_{k-s}) \land C_1(x_1, x_{k-s})))))),
$$
$$
Q_2(x_1) = \exists x_2 \mbox{ } (D^{=}_{a_2}(x_1, x_2) \land (\exists x_3 \mbox{ } (D^{=}_{2^{k-s-1}}(x_2, x_3) 
\land \ldots \land (\exists x_{k-s} \mbox{ } (D^{=}_{2^2}(x_{k-s-1}, x_{k-s}) \land C_2(x_1, x_{k-s})))))),
$$
\begin{multline*}
C_1(x_1, x_{k-s}) = ((x_{k-s} \neq x_1) \land (\exists x_{k - s + 1} \exists x_{k - s + 2} \mbox{ } ((x_{k - s + 1} \neq x_1) \land (x_{k - s + 2} \neq x_1) \land \\ \land
(\exists y_1 \ldots \exists y_{s-3} \mbox{ } (N(x_{k-s}, x_{k-s+1}, x_{k-s+2}, y_1, \ldots, y_{s-3}) \land (y_1 \neq x_1) \land \ldots \land (y_{s-3} \neq x_1))) \land \\
\land (\exists y_1 \ldots \exists y_{s-2} \mbox{ } (N(x_{k-s}, x_{k-s+1}, y_1, \ldots, y_{s-2}) \land (y_1 \neq x_1) \land (y_1 \neq x_{k-s+2}) \land \ldots \land \\
\land (y_{s-2} \neq x_1) \land (y_{s-2} \neq x_{k-s+2})))))),
\end{multline*}
\begin{multline*}
C_2(x_1, x_{k-s}) = ((x_{k-s} \neq x_1) \land (\exists x_{k - s + 1} \exists x_{k - s + 2} \mbox{ } ((x_{k - s + 1} \neq x_1) \land (x_{k - s + 2} \neq x_1) \land R_1(x_1, x_{k-s}, x_{k-s+1}) \land \\ \land R_1(x_1, x_{k-s}, x_{k-s+2}) \land R_1(x_1, x_{k-s+1}, x_{k-s+2}) \land R_2(x_{k-s}, x_{k-s+1}, x_{k-s+2})))),
\end{multline*}
$$
R_1(x, y_1, y_2) = \exists y_3 \ldots \exists y_s \mbox{ } (N(y_1, y_2, y_3, \ldots, y_s) \land (y_3 \neq x) \land \ldots \land (y_s \neq x)),
$$
$$
R_2(y_1, y_2, y_3) = \neg (\exists y_4 \ldots \exists y_s \mbox{ } N(y_1, y_2, y_3, \ldots,  y_s))
$$
for $s > 3$ and $R_2(y_1, y_2, y_3) = \neg N(y_1, y_2, y_3)$ for $s = 3$.

The formula $D^{=}_i(x_1, x_2)$ was defined in the proof of Theorem 6 and expresses the property that the distance between the vertices $x_1$ and $x_2$ equals $i$. The quantifier depth of the formula $D^{=}_i(x_1, x_2)$ equals $\lceil \log_2 i \rceil + s - 2$.

Let $H_1$ be a hypergraph with $V(H_1) = \{x_1^1, \ldots, x_{2(s-1)}^1\}$ and $E(H_1) = \{\{x_1^1, \ldots, x_s^1\}$, $\{x_s^1, \ldots, x_{2s-2}^1, x_1^1\}\}$. 
Let $H_2$ be a hypergraph with $V(H_2) = \{x_1^2, \ldots, x_{3(s-1)}^2\}$ and 
$E(H_2) = \{\{x_1^2, \ldots, x_s^2\}, \{x_s^2, \ldots, x_{2s-1}^2\}$, $\{x_{2s-1}^2, \ldots, x_{3s-3}^2, x_1^2\}\}$.
Let $H$ be a hypergraph obtained from the disjoint union of the hypergraphs $H_1$ and $H_2$ by adding a vertex $x$ and two non-intersecting paths with lengths $a_1 + 2^{k-s} - 4$ and $a_2 + 2^{k-s} - 4$ which connect the vertex $x$ and the vertices $x_1^1$ and $x_1^2$ respectively. 

Let $\tilde \Omega_n$ be the set of all hypergraphs $\mathcal G$ from $\Omega_n$ such that for any hypergraph $G$ with $\rho^{\max}(G) > \frac{1}{\alpha}$ and 
$v(G) \le 2^k s$, there is no copy of $G$ in $\mathcal G$.

Let us prove that if the hypergraph $\mathcal G \in \tilde \Omega_n$ satisfies $L$ then it contains a copy of the hypergraph $H$. 
Suppose that the hypergraph $\mathcal G \in \tilde \Omega_n$ satisfies $L$. Then there exist vertices $x_1, x_2^1, x_2^2 \in V(\mathcal G)$ such that $d_{\mathcal G}(x_1, x_2^1) \le a_1 + 2^{k-s-1} + \ldots + 2^2 = a_1 + 2^{k-s} - 2^2$, $d_{\mathcal G}(x_1, x_2^2) \le a_2 + 2^{k-s-1} + \ldots + 2^2 = a_2 + 2^{k-s} - 2^2$ and the predicates $C_1(x_1, x_2^1)$ and $C_2(x_1, x_2^2)$ are true. As $\mathcal G \in \tilde \Omega_n$, the truth of the predicate $C_1(x_1, x_2^1)$ implies that there exists a subhypergraph $\tilde H_1$ in $\mathcal G$ which is isomorphic to the hypergraph $H_1$, contains the vertex $x_2^1$ and does not contain
the vertex $x_1$. The truth of the predicate $C_2(x_1, x_2^2)$ implies that there exists a subhypergraph $\tilde H_2$ in $\mathcal G$ which is isomorphic to the hypergraph $H_2$, contains the vertex $x_2^2$ and does not contain
the vertex $x_1$. Moreover, we have $V(\tilde H_1) \cap V(\tilde H_2) = \varnothing$, since otherwise the density of the hypergraph $\tilde H_1 \cup \tilde H_2$ is greater than $1/\alpha$, which contradicts $\mathcal G \in \tilde \Omega_n$. Let $P_1$ be the shortest path which connects the vertex $x_1$ with some vertex from $\tilde H_1$ and $P_2$ be the shortest path which connects the vertex $x_1$ with some vertex from $\tilde H_2$. Let $X_0 = \tilde H_1 \cup \tilde H_2$, $X_1 = \tilde H_1 \cup \tilde H_2 \cup P_1$ and $X_2 = \tilde H_1 \cup \tilde H_2 \cup P_1 \cup P_2$. 
Then $e(X_1, X_0) \le a_1 + 2^{k-s} - 4$ and $v(X_1, X_0) \le e(X_1, X_0)(s-1)$. If $P_2 \subset \tilde H_1 \cup \tilde H_2 \cup P_1$, then $v(X_1, X_0) \le e(X_1, X_0)(s-1) - 1$.
If $P_2 \not \subset \tilde H_1 \cup \tilde H_2 \cup P_1$, then $e(X_2, X_1) \le a_2 + 2^{k-s} - 4$ and $v(X_2, X_1) \le e(X_2, X_1)(s-1) - 1$. Thus $|E(X_2)| = e(X_1, X_0) + e(X_2, X_1) + 5 \le 2^{k-s+1} + a$ and $|V(X_2)| = v(X_1, X_0) + v(X_2, X_1) + 5(s-1) \le |E(X_2)| (s - 1) - 1$. Therefore, 
$$1/\rho(X_2) \le s - 1 - \frac{1}{|E(X_2)|} \le s - 1 - \frac{1}{2^{k-s+1} + a} = \alpha.$$ 
Equalities hold if and only if $e(X_i, X_{i-1}) = a_i + 2^{k-s} - 4$ for any $i \in \{1, 2\}$ and $v(X_1, X_0) = (a_1 + 2^{k-s} - 4) (s-1)$, $v(X_2, X_1) = (a_2 + 2^{k-s} - 4) (s-1) - 1$. 
Therefore, by the definition of $\tilde \Omega_n$ these equalities hold and the hypergraph $X_2$ is isomorphic to $H$.

Conversely, if the hypergraph $\mathcal G \in \tilde \Omega_n$ contains a copy of the hypergraph $H$ then $\mathcal G$ satisfies $L$. 
Since $\lim_{n \to \infty} \Pr[G^s(n, n^{-\alpha}) \in \tilde \Omega_{n}] = 1$, it follows that $\lim_{n \to \infty} \Pr[G^s(n, n^{-\alpha}) \models L] = \lim_{n \to \infty} \Pr[G^s(n, n^{-\alpha}) \models L_H] \in (0, 1)$. Therefore, $G^s(n, n^{-\alpha})$ does not obey the zero-one $k$-law. 

Finally, consider the case $k = s + 1$. Then $a = 1$ and $\alpha = s - 1 - \frac{1}{5}$. 
Let $L$ be a first-order property which is expressed by the formula $\exists x_1 \mbox{ } (Q_1(x_1) \land Q_2(x_1))$, where
\begin{multline*}
Q_1(x_1) = \exists x_2 \exists x_3 \mbox{ } ((\exists y_1 \ldots \exists y_{s-3} \mbox{ } N(x_1, x_2, x_3, y_1, \ldots, y_{s-3})) \land \\
\land (\exists y_1 \ldots \exists y_{s-2} \mbox{ } (N(x_1, x_2, y_1, \ldots, y_{s-2}) \land (y_1 \neq x_3) \land \ldots \land (y_{s-2} \neq x_3)))),
\end{multline*}
$$
Q_2(x_1) = \exists x_2 \exists x_3 \mbox{ } (T_2(x_1, x_2) \land T_2(x_1, x_3) \land T_2(x_2, x_3) \land (\neg T_3(x_1, x_2, x_3))),
$$
$$
T_i(y_1, \ldots, y_i) = \exists y_{i+1} \ldots \exists y_s \mbox{ } N(y_1, \ldots, y_i, y_{i+1}, \ldots, y_s)
$$
for $i < s$ and $T_s(y_1, \ldots, y_s) = N(y_1, \ldots, y_s)$.

Let $H$ be a hypergraph with $V(H) = \{x_1, x_2^1, \ldots, x_{2(s-1)}^1, x_2^2, \ldots, x_{3(s-1)}^2\}$ and 
$$E(H) = \{\{x_1,x_2^1 \ldots, x_s^1\}, \{x_s^1, \ldots, x_{2s-2}^1, x_1\}, \{x_1, x_2^2 \ldots, x_s^2\}, \{x_s^2, \ldots, x_{2s-1}^2\}, \{x_{2s-1}^2, \ldots, x_{3s-3}^2, x_1\}\}.$$ 

Let $\tilde \Omega_n$ be the set of all hypergraphs $\mathcal G$ from $\Omega_n$ such that for any hypergraph $G$ with $\rho^{\max}(G) > \frac{1}{\alpha}$ and 
$v(G) \le 2^k s$, there is no copy of $G$ in $\mathcal G$. 

Let us show that if the hypergraph $\mathcal G \in \tilde \Omega_n$ satisfies $L$ then it contains a copy of the hypergraph $H$. 
Suppose that the hypergraph $\mathcal G \in \tilde \Omega_n$ satisfies $L$. Then there exists a vertex $x_1$ and there exist subhypergraphs $\tilde H_1$ and $\tilde H_2$ in $\mathcal G$ which are isomorphic to the hypergraphs $H_1$ and $H_2$ respectively and $x_1 \in V(\tilde H_1) \cap V(\tilde H_2)$. If $V(\tilde H_1) \cap V(\tilde H_2) \neq \{x_1\}$, then the density of the hypergraph $\tilde H_1 \cup \tilde H_2$ is greater than $1/\alpha$, which contradicts $\mathcal G \in \tilde \Omega_n$. Therefore, the hypergraph $\tilde H_1 \cup \tilde H_2$ is isomorphic to $H$. Obviously, if the hypergraph $\mathcal G \in \tilde \Omega_n$ contains a copy of the hypergraph $H$, then it satisfies $L$.   
Since $\lim_{n \to \infty} \Pr[G^s(n, n^{-\alpha}) \in \tilde \Omega_{n}] = 1$, we obtain that $\lim_{n \to \infty} \Pr[G^s(n, n^{-\alpha}) \models L] = \lim_{n \to \infty} \Pr[G^s(n, n^{-\alpha}) \models L_H] \in (0, 1)$.
\\

\textbf{Acknowledgements}
\\

This work has been supported by the grant NSh-775.2022.1.1. 

\renewcommand{\refname}{References}


\begin{thebibliography}{50}

\bibitem{shelah_spencer}
S. Shelah, J.H. Spencer, Zero-one laws for sparse random graphs, J. Amer. Math. Soc. 1: 97-115, 1988.

\bibitem{inf_spectra} J.H. Spencer, Ininite spectra in the first order theory of graphs, Combinatorica 10(1): 95-102, 1990.

\bibitem{zhuk1}
M.E. Zhukovskii, Zero-one k-law, Discrete Mathematics, 2012, 312: 1670-1688.

\bibitem{zhuk2}
M.E. Zhukovskii, On the zero-one k-law extensions, European J. of Combinatorics, 60 (2017): 66–81.

\bibitem{zhuk_spencer}
Spencer, J., Zhukovskii, M. E. Bounded quantifier depth spectra for random graphs. Discrete Mathematics, 339(6), 1651-1664, 2016.

\bibitem{zhuk4}
M.E. Zhukovskii, “On infinite spectra of first order properties of random graphs”, Moscow Journal of Combinatorics and Number Theory, 2016, Vol. 6, No. 4, 73–102.

\bibitem{zhuk5}
M.E. Zhukovskii, “Logical laws for short existential monadic second order sentences about graphs”, Journal of Mathematical Logic, Vol 2, No 2, 2050007 (2020).

\bibitem{zhuk_popova}
S.N. Popova, M.E. Zhukovskii, “Existential monadic second order logic of undirected graphs: a disproof of the Le Bars conjecture”, Annals of Pure and Applied Logic, 170 (2019) 505–514.

\bibitem{strange_logic} J. H. Spencer, The Strange Logic of Random Graphs, Number 22 in Algorithms and Combinatorics,
Springer-Verlag, Berlin, 2001.

\bibitem{raigor}
M.E. Zhukovskii, A.M. Raigorodskii, Random graphs: models and asymptotic characteristics,
Russian Mathematical Surveys (2015), 70(1): 33.

\bibitem{ebbing}
H.D. Ebbinghaus and J. Flum. Finite model theory. Perspectives in Mathematical Logic. Springer-Verlag, Berlin, second edition, 1999. 

\bibitem{zero_one}
A.D. Matushkin, S.N. Popova, Strictly balanced uniform hypergraphs and generalizations of Zero-One Law,
Discrete Mathematics, 345(6), 2022, 112835.

\bibitem{infspectrum} S.N. Popova, Infinite spectra of first-order properties for random hypergraphs, 
Problems of Information Transmission, 54(3): 281-289, 2018.

\bibitem{spectrum}
S.N. Popova, Spectrum for first-order properties of random hypergraphs, arXiv:1908.01074.

\bibitem{limitspectrum}
S.N. Popova, Limit points of spectra for first-order properties of random hypergraphs, Discrete Applied Mathematics, 293: 134-142, 2021. 

\bibitem{vantsyan}
A.G. Vantsyan, The evolution of random uniform hypergraphs, Probabilistic problems
in discrete mathematics, 126-131, Moscov. Inst. Electron. Machinostroenya, 1987.

\bibitem{JLR}
S. Janson, T. \L uczak, A. Rucinski, Random Graphs, New York, Wiley, 2000.

\bibitem{alon_spencer} 
N. Alon and J.H. Spencer, The Probabilistic Method, Wiley-Interscience, New York, 2000.

\bibitem{zhuk_ext} M.E. Zhukovskii, Estimation of the number of maximal extensions in a random graph, Discrete Mathematics and Applications, 22:1, 55–90, 2012.


\end{thebibliography}
\end{document}